\documentclass{endm}
\usepackage{endmmacro}
\usepackage{graphicx}
\usepackage{enumerate}

\def\diam{{\rm diam}}
\def\span{{\rm span}}
\def\rn{{\rm rn}}

\def\l{\lambda}

\def\ve{\varepsilon}

\begin{document}

\begin{verbatim}\end{verbatim}\vspace{2.5cm}

\begin{frontmatter}

\title{Further results on the radio number of trees}

\author{Devsi Bantva\thanksref{myemail}}
\address{Department of Mathematics\\ Lukhdhirji Engineering College, Morvi - 363 642, Gujarat, India}

\thanks[myemail]{Email:\href{mailto:devsi.bantva@gmail.com} {\texttt{\normalshape devsi.bantva@gmail.com}} (Devsi Bantva)}

\begin{abstract}
Let $G$ be a finite, connected, undirected graph with diameter $\diam(G)$ and $d(u,v)$ denote the distance between $u$ and $v$ in $G$. A radio labeling of a graph $G$ is a mapping $f: V(G) \rightarrow \{0,1,2,...\}$ such that $|f(u)-f(v)| \geq \diam(G) + 1 - d(u,v)$ for every pair of distinct vertices $u, v$ of $G$. The radio number of $G$, denoted by $\rn(G)$, is the smallest integer $k$ such that $G$ has a radio labeling $f$ with $\max\{f(v) : v \in V(G)\} = k$. In this paper, we determine the radio number for three families of trees obtained by taking graph operation on a given tree or a family of trees.

\end{abstract}


\end{frontmatter}

\section{Introduction}\label{intro}

A number of graph labelings have a root in channel assignment problem. In the channel assignment problem, we seek to assign channels to transmitters such that it satisfies all interference constraints. This well explored problem is also studied using graph coloring. In a graph, a set of transmitters is represented by vertices of a graph; two vertices are adjacent if transmitters are \emph{very close} and at distance two apart if transmitters are \emph{close} in a network. Notice that two transmitters are classified as \emph{very close} if the interference is unavoidable and \emph{close} if the interference is avoidable between them. Motivated by this problem Griggs and Yeh \cite{Griggs} introduced the following $L(2,1)$-labeling problem: An \emph{$L(2,1)$-labeling} (or \emph{distance-two labeling}) of a graph $G=(V(G),E(G))$ is a function $f$ from the vertex set $V(G)$ to the set of non-negative integers such that $|f(u)-f(v)|\geq2$ if $d(u,v)=1$ and $|f(u)-f(v)|\geq1$ if $d(u,v)=2$. The \emph{span} of $f$ is defined as $\max\{|f(u)-f(v)| : u, v \in V(G)\}$, and the minimum span over all $L(2,1)$-labelings of $G$ is called the \emph{$\lambda$-number} of $G$, denoted by $\l(G)$. Observe that $L(2,1)$-labeling deal with two level interference only. The $L(2,1)$-labeling and other distance-two labeling problems have been studied by many researchers in the past two decades; see \cite{Calamoneri} and \cite{Yeh1}.

In 2005, Chartrand \emph{et al.} \cite{Chartrand1} introduced the concept of radio labeling and put the level of interference at largest possible-the diameter of graph. Denote by $\diam(G)$ the \emph{diameter} of $G$, that is, the maximum distance among all pairs of vertices in $G$.

\begin{definition}
A \emph{radio labeling} of a graph $G$ is a mapping $f: V(G) \rightarrow \{0, 1, 2, \ldots\}$ such that for every pair of distinct vertices $u, v$ of $G$,
$$
d(u,v) + |f(u)-f(v)| \geq \diam(G) + 1.
$$
The integer $f(u)$ is called the \emph{label} of $u$ under $f$, and the \emph{span} of $f$ is defined as $\span(f) = \max \{|f(u)-f(v)|: u, v \in V(G)\}$. The \emph{radio number} of $G$ is defined as
$$
\rn(G) := \min_{f} \span(f)
$$
with minimum taken over all radio labelings $f$ of $G$. A radio labeling $f$ of $G$ is \emph{optimal} if $\span(f) = \rn(G)$.
\end{definition}

Note that any optimal radio labeling must assign $0$ to some vertex and also in the case when $\diam(G) = 2$ we have $\rn(G) = \l(G)$. Observe that any radio labeling should assign different labels to distinct vertices. In fact, a radio labeling induces an ordering $u_{0}$, $u_{1}$ ,...,$u_{p-1}$ ($p$ = $|V(G)|$) of vertices such that $0 = f(u_{0}) < f(u_{1}) < ... < f(u_{p-1}) = \span(f)$.

The radio number of graphs is studied by limited group of authors. The readers are advised to refer the following papers for the radio number of listed graph families; \cite{Chartrand1,Chartrand2,Liu,Zhang} for paths and cycles, \cite{Daphne2,Daphne3} for square of paths and cycles, \cite{Benson} for all graphs of order $n$ and diameter $n-2$, \cite{Cada} for distance graphs, \cite{Bantva1,Bantva2,Daphne1} about radio number of trees, \cite{Li} for complete $m$-ary trees, \cite{Tuza} for level-wise regular trees, \cite{Vaidya1} for total graph of paths, \cite{Vaidya2} for strong product $P_{2} \boxtimes P_{n}$, \cite{Vaidya3} for linear cacti. One can also refer to the survey article \cite{Chartrand} for detail on the radio number of graphs.

The results presented in this paper are in continuation of our previous work in \cite{Bantva1,Bantva2}. If a tree or a family of trees satisfies condition of \cite[Theorem 3.2]{Bantva2} then we present three families of trees obtained by taking graph operation on a given tree or a family of trees such that it also satisfies condition of \cite[Theorem 3.2]{Bantva2}. We relate the radio number of it with the radio number of the given tree or family of trees. Proofs of our results will be given in the full version of this paper.

\section{Preliminaries and earlier work on \rn(T)}

A \emph{tree} $T$ is a connected graph that contains no cycle. For a tree $T$, denote \emph{vertex set} and \emph{edge set} by $V(T)$ and $E(T)$. The \emph{order} of a tree $T$ is the number of vertices in it. The \emph{distance} $d(u,v)$ between two vertices $u$ and $v$ is the length of a shortest path connecting them. The \emph{diameter} of a graph $G$ is max\{$d(u,v)$ : $u, v \in V(G)$\}. For a vertex $v \in V(G)$, the \emph{neighborhood of $v$} denoted by $N(v)$, is the set of vertices adjacent to $v$. Terms and notations not defined here are used in the sense of \cite{West}.

The first result on the radio number of trees was given by Chartrand \emph{et al.} in \cite{Chartrand1,Chartrand2}. They gave an upper bound for the radio number of paths and trees. Later, Liu and Zhu gave the exact radio number of paths in \cite{Liu}. The lower bound for the radio number of trees and different necessary and sufficient condition to achieve the lower bound is given by Liu \cite{Daphne1} and Bantva \emph{et al.} \cite{Bantva2}. The radio number for particular trees is determined by many authors; see \cite{Bantva1,Benson,Tuza,Li}. In spite of these efforts, the problem of determining the exact value of the radio number for trees is still open. However, for present work, our main concern is with \cite{Daphne1} and \cite{Bantva1,Bantva2} and hence we present the results of it with necessary terms and notations.

In \cite{Daphne1}, author gave a lower bound for the radio number of trees and presented necessary and sufficient condition to achieve this lower bound. She also presented a class of trees, namely spiders, achieving this lower bound in \cite{Daphne1}. She defined several terms and notations in \cite{Daphne1} to give this lower bound. First we present these terms and notations which are as follows. In \cite{Daphne1}, author viewed $T$ as  rooted at a single vertex $w$ and defined the \emph{level function} on $V(T)$ from fix root $w$ by $L_{w}(u)$ = $d(w,u)$ for any $u \in V(T)$. For any two vertices $u$ and $v$, if $u$ is on the ($w$, $v$)-path, then $u$ is an \emph{ancestor} of $v$, and $v$ is \emph{descendent} of $u$. If $u$ be a neighbor of $w$ then the subtree induced by $u$ together with all the descendent of $u$ is called a \emph{branch} at $u$. Two branches are called \emph{different} if they are induced by two different vertices adjacent to $w$. The \emph{weight of $T$} from $v \in V(T)$ is defined as $w_{T}(v)$ = $\sum_{v \in V(T)} d(u,v)$ and the \emph{weight of $T$} as $w(T)$ = min\{$w_{T}(v)$ : $v \in V(T)$\}. A vertex $v \in V(T)$ is a \emph{weight center} of $T$ if $w_{T}(v)$ = $w(T)$. We denote the set of weight centers of $T$ by $W(T)$. In \cite{Daphne1}, it is proved that every tree $T$ has either one or two weight centers, and $T$ has two weight centers, say, $W(T)$ = \{$w$, $w^{'}$\}, if and only if $w$ and $w^{'}$ are adjacent and $T-ww^{'}$ consists of two equal-sized components. Using these terms and notations Liu presented the following result in \cite{Daphne1}.
\begin{theorem}\cite{Daphne1} Let $T$ be an $m$-vertex tree with diameter $d$. Then
\begin{equation}
\rn(T) \geq (m-1)(d+1)+1-2w(T).
\end{equation}
Moreover, the equality holds if and only if for every weight center $w^{*}$, there exist a radio labeling $f$ with $f(u_{0})$ = 0 $<$ $f(u_{1})$ $<$ ... $<$ $f(u_{m-1})$, where all the following hold (for all $0 \leq i \leq m-2$); \\
(1) $u_{i}$ and $u_{i+1}$ are in different branches (unless one of them is $w^{*}$); \\
(2) \{$u_{0}$, $u_{m-1}$\} = \{$w^{*}$, $v$\}, where $v$ is some vertex with $L_{w^{*}}(v)$ = 1; \\
(3) $f(u_{i+1})$ = $f(u_{i})+d+1-L_{w^{*}}(u_{i})-L_{w^{*}}(u_{i+1})$.
\end{theorem}

Bantva \emph{et al.}\cite{Bantva2} modified the lower bound given by Liu and gave more useful necessary and sufficient condition. They viewed $T$ as rooted at its weight center $W(T)$: if $W(T)$ = \{w\}, then $T$ is rooted at $w$; if $W(T)$ = \{$w$, $w^{'}$\} (where $w$ and $w^{'}$ are adjacent), then $T$ is rooted at $w$ and $w^{'}$ in the sense that both $w$ and $w^{'}$ are at level 0. They called two branches are \emph{different} if they are at two vertices adjacent to the same weight center (which is same as in \cite{Daphne1}), and \emph{opposite} if they are at two vertices adjacent to different weight centers. The later case occurs only when $T$ has two weight centers. They defined the \emph{level of $u$} in $T$ as
\begin{equation}
L(u) := \mbox{min}\{d(u,x) : x \in W(T)\}, u \in V(T)
\end{equation}
and the \emph{total level of $T$} as
\begin{equation}
L(T) := \displaystyle\sum_{u \in V(T)} L(u).
\end{equation}
Define
\begin{eqnarray*}
\ve(T) & = &\left\{
\begin{array}{l}
\begin{tabular}{lll}
1, & \mbox{ if $T$ has only one weight center}, \\
0, & \mbox{ if $T$ has two (adjacent) weight centers}.
\end{tabular}
\end{array}
\right.
\end{eqnarray*}
Using these terms and notations, Bantva \emph{et al.} presented the following results in \cite{Bantva2}.
\begin{lemma}\cite{Bantva2}\label{thm:lb}
Let $T$ be a tree with order $p$ and diameter $d \ge 2$. Denote $\ve = \ve(T)$. Then
\begin{equation}\label{eq:lb}
\rn(T) \ge (p-1)(d+\ve) - 2 L(T) + \ve.
\end{equation}
\end{lemma}

\begin{theorem}\cite{Bantva2}\label{thm:ub}
Let $T$ be a tree with order $p$ and diameter $d \ge 2$. Denote $\ve = \ve(T)$. Then
\begin{equation}\label{eq:ub}
\rn(T) = (p-1)(d+\ve) - 2 L(T) + \ve
\end{equation}
holds if and only if there exists a linear order $u_{0}, u_{1}, \ldots, u_{p-1}$ of the vertices of $T$ such that
\begin{enumerate}[(a)]
\item $u_{0} = w$ and $u_{p-1} \in N(w)$ when $W(T) = \{w\}$, and $\{u_{0}, u_{p-1}\} = \{w, w^{'}\}$ when $W(T) = \{w, w^{'}\}$;
\item the distance $d(u_{i}, u_{j})$ between $u_{i}$ and $u_{j}$ in $T$ satisfies ($0 \leq i < j \leq p-1$)
\end{enumerate}
\begin{equation}\label{eq:dij}
d(u_{i},u_{j}) \geq \displaystyle\sum_{t = i}^{j-1} (L(u_{t})+L(u_{t+1})) - (j - i)(d+\ve) + (d+1).
\end{equation}
Moreover, under this condition the mapping $f$ defined by
\begin{equation}\label{eq:f0}
f(u_{0}) = 0
\end{equation}
\begin{equation}\label{eq:f}
f(u_{i+1}) = f(u_{i}) - L(u_{i+1}) - L(u_{i}) + (d + \ve),\;\, 0 \leq i \leq p-2
\end{equation}
is an optimal radio labeling of $T$.
\end{theorem}

Using Theorem \ref{thm:ub}, Bantva \emph{et al.} determined the radio number for banana trees, firecracker trees and a special class of caterpillars in \cite{Bantva2}. Bantva \emph{et al.} also noticed in \cite{Bantva2} that the radio number of paths and complete $m$-ary tree can also be determined using Theorem \ref{thm:ub}.

\section{Main results}

In this section, we present some more results on radio number of trees using Theorem \ref{thm:ub}. For this purpose we continue to use the terminology and notation defined in the previous section.

Now we consider a tree $T$ of order $n_{0}$ and diameter $d_{0}$ with weight center $w_{0}$. In case of families of trees, we consider trees $T_{i}$ ($1 \leq i \leq k$) of order $n_{i}$ and diameter $d_{i}$ with weight center $w_{i}$. If $T_{*}$ is any tree obtained by taking graph operation on $T$ or family of trees $T_{i}$ then we take $|T_{*}|$ = $n$ and diam($T_{*}$) = $d$. A $k$-star is a tree consisting of $k$ leaves and another vertex joined to all leaves by edges. A $k$-double star is a tree which is formed by joining $k$ edges to each of the two vertices of $K_{2}$. We define $T_{w_{k}}$ to be the tree obtained by identifying weight center $w_{i}$ of trees $T_{i}$, $1 \leq i \leq k$ with a single vertex $w$. Note that the weight center of $T_{w_{k}}$ is $w$ and $|T_{w_{k}}|$ = $kn_{i}-k+1$. We define $T_{S_{k}}$ to be the tree obtained by taking $k$ copies of a tree $T$ and identifying a weight center $w$ of each of them with each leaf of $k$-star. The weight center of $T_{S_{k}}$ is $w$ and $|T_{S_{k}}|$ = $kn_{0}+1$. We define $T_{D_{k}}$ to be the tree obtained by taking $2k$ copies of tree $T$ and identifying a weight center $w$ of each of them with each leaf of $k$-double star. Note that $T_{D_{k}}$ has two adjacent weight centers and $|T_{D_{k}}|$ = $2(kn_{0}+1)$.

\begin{theorem}\label{thm:Twk} If Theorem \ref{thm:ub} hold for $T_{i}, 1 \leq i \leq k$ then Theorem \ref{thm:ub} holds for $T_{w_{k}}$ and \begin{equation}
\rn(T_{w_{k}}) = \displaystyle\sum_{i=1}^{k} [\rn(T_{i})+(n_{i}-1)(d-d_{i})]-k+1.
\end{equation}
\end{theorem}

\begin{theorem}\label{thm:Tsk} If Theorem \ref{thm:ub} holds for $T$ then Theorem \ref{thm:ub} holds for $T_{S_{k}}$ and
\begin{equation}
\rn(T_{S_{k}}) = k[\rn(T)+n_{0}(d-d_{0}-2)+d_{0}]+1.
\end{equation}
\end{theorem}

\begin{theorem}\label{thm:Tpk} If Theorem \ref{thm:ub} holds for $T$ then Theorem \ref{thm:ub} holds for $T_{D_{k}}$ and
\begin{equation}
\rn(T_{D_{k}}) = 2k[\rn(T)+n_{0}(d-d_{0}-3)+d_{0}]+d.
\end{equation}
\end{theorem}
Proofs and illustrations of Theorem \ref{thm:Twk} to \ref{thm:Tpk} will be provided in full version of this paper.

\section*{Acknowledgements}
The author is grateful to the reviewers for their corrections and suggestions. The author also thank Samir Vaidya and Sanming Zhou for earlier discussion on this topic.

\end{document}